\documentclass[preprint, 12pt]{article}

\pdfoutput=1

\usepackage{amsmath,amsthm,amsfonts,amssymb,amscd}
\usepackage{graphicx,bm,color}
\usepackage{algorithm}
\usepackage[algo2e]{algorithm2e}
\usepackage{amsfonts}
\usepackage{algpseudocode}
\usepackage{stmaryrd}
\usepackage{listings}
\usepackage{lineno}
\usepackage{mathptmx}
\usepackage{cleveref}
\usepackage{tikz}%For some graphics stuff if required
\usepackage[skins,theorems]{tcolorbox}
\usepackage{authblk}
\usepackage{tkz-graph}
\usepackage[latin1]{inputenc}
\usetikzlibrary{shapes,arrows}
\usepackage{smartdiagram}
\usetikzlibrary{arrows,calc, decorations.markings, intersections}
\usepackage{caption}
\usepackage{cancel}
\usepackage{mathtools}
\usepackage{url}
\usepackage{float}%To place a float exactly at a point (using H instead of !ht)
 
\definecolor{gray}{gray}{0.6}

\theoremstyle{plain}
\theoremstyle{remark}

\newtheorem*{remark*}{Remark}

\newcommand{\vertiii}[1]{{\left\vert\kern-0.25ex\left\vert\kern-0.25ex\left\vert #1 
		\right\vert\kern-0.25ex\right\vert\kern-0.25ex\right\vert}}

\newcommand\norm[1]{\left\lVert#1\right\rVert}

\begin{document}

\title{Recovery of the Interface Velocity for the Incompressible Flow in Enhanced Velocity \\ Mixed Finite Element Method}

%\author{Yerlan Amanbek  \and Gurpreet Singh \and Gergina Pencheva \and Mary F. Wheeler
%}

\author[1,2]{Yerlan Amanbek}
\author[1]{Gurpreet Singh}
\author[1]{Mary F. Wheeler}
\affil[1]{Institute for Computational Engineering and Sciences, University of Texas at Austin}
\affil[2]{Nazarbayev University}
\affil[ ]{\textit{yerlan.amanbek@nu.edu.kz, \{gurpreet, mfw\}@ices.utexas.edu}}

%\author[Affil1]{First Author}
%\author[Affil2]{Second Author}
%%\author[Affil1]{Third/Corresponding Author \corref{cor1}}
%\address[Affil1]{Affiliation address 1}
%\address[Affil2]{Affiliation address 2}
%% Replace capitalized text with the appropriate information (use standard capitalization rules for your text, not all capitals.
%\cortext[cor1]{Corresponding Author: AUTHOR'S NAME, TYPE AUTHOR'S POSTAL ADDRESS; Email, TYPE CORRESPONDING AUTHOR'S EMAIL ADDRESS; Phone, TYPE CORRESPONDING AUTHOR'S PHONE NUMBER}

\date{\today}

\maketitle              % typeset the header of the contribution
\begin{abstract}
	The velocity, coupling term in the flow and transport problems, is important in the accurate numerical simulation or in the \textit{posteriori} error analysis for adaptive mesh refinement. We consider Enhanced Velocity Mixed Finite Element Method for the incompressible Darcy flow.
	In this paper, our aim to study the improvement of velocity at interface to achieve the better approximation of velocity between subdomains. 
	We propose the reconstruction of velocity at interface by using the post-processed pressure.
	Numerical results at the interface show improvement on convergence rate.
\end{abstract}

\begin{keyword}
	domain decomposition, velocity improvement, enhanced velocity, mixed fem.
\end{keyword}

%------------------------------------------------------------------------------
\section{Introduction}
The numerical reservoir simulations have been utilized in many subsurface applications such as groundwater remediation, reservoir well evaluation, and contaminate transport problems. For such applications, it is common to deal with the flow and transport problem. The main component or coupling term of the flow and transport systems is the velocity and its accuracy the mostly achieved by employing classical mixed finite element system. Due to the heterogeneity of porous media multiphysics problems could be categorized systematically in which one physical phenomena influences within a subdomain and another physical phenomena dominates within another subdomain. Such solutions are coupled through continuity of normal flux at interface, shared region between differently discretized subdomains. To deal with these problems there are the well-known methods such as Multiscale Mortar and Enhanced Velocity schemes that are established in various applications. Recently, a novel adaptive method was studied in subsurface applications \cite{amanbek2017adaptive,singh2017adaptive,amanbek2018priori,amanbek2018new} using Enhanced Velocity scheme. The main idea is here to utilize the EVMFEM as domain decomposition method to couple different discretized subdomains with more accurate upscaled subsurface parameters.

In the simulation of flow with adaptivity, the results obtained in  \cite{gerritsen2008integration} suggest that pressure values could be interpolated using neighboring elements values to approximate auxiliary pressure values within provided elements. Selection of interpolants is based on convex combinations of vertical and horizontal oriented pressure values. 
In related reference \cite{arbogast2014posteriori}, it was studied that the interface error of solution between subdomains for different numerical methods including Mortar Multiscale MixedFEM which provided better approximation for second-order elliptic problems. One of reasons is the iterative procedure in the mortar scheme that is a key in coupling two subdomains physics. According to author in \cite{arbogast2014posteriori} mortar scheme is general method in coupling for practical multiphysics problems. On the other hand, the efficient Enhanced Velocity scheme has not been investigated from the point of view of the improvement solution including velocity at interface in the previous studies.

The challenge here is to construct the velocity approximation of EVMFEM and specifically at interface to have a better velocity between subdomains that leads accurate approximation in the flow and transport problems. In \cite{wheeler2002enhanced}, a \textit{priori} error analysis states that the global error is 
\begin{align}
\norm{\mathbf{u}-\mathbf{u}_h}_{\Omega} \le C \left(\norm{p}_{1, \infty, \Omega^*} + \norm{\mathbf{u}}_{1, \Omega} h^{1/2}\right) h^{1/2} \label{eq:fluxerrordomain}
\end{align}
and away from the interface $\Gamma$ the velocity error convergence rate is better, since
\begin{align}
\norm{\mathbf{u}-\mathbf{u}_h}_{\Omega^{'}} \le C_{\varepsilon} \left(\norm{p}_{1, \infty, \Omega^*} + \norm{\mathbf{u}}_{1, \Omega} \right) h^{r-\varepsilon}  \label{eq:fluxinsidedomain}
\end{align}
where $\varepsilon > 0$, $r=1$ if $d=2$ and $r=5/6$ if $d=3$, and $\Omega^{'}_i$ is compactly contained in $\Omega_i$, $\Omega^{'} = \bigcup_{i=1}^{N_b} \Omega^{'}_i$. This implies that the discrete velocity should be approximated more precise near interface region $\Omega^*$. On the question of pressure approximation, the convergence rate of pressure approximation is $\mathcal{O}(h^{1})$, if $d=2$, and $\mathcal{O}(h^{5/6})$, if $d=3$ \cite{wheeler2002enhanced,thomas2011enhanced}. If one compare the error of velocity (\ref{eq:fluxerrordomain}) and pressure approximation these results indicate that the velocity convergence rate is not strong as pressure in $\Omega$.  Similar a \textit{priori} error result was shown in \cite{amanbek2018priori} for transient problems. Nevertheless, there are still problems including the velocity approximation at the interface to be addressed.

In this paper, we introduce the way to improve velocity accuracy at interface in the Enhanced Velocity MFEM for incompressible flow using the post-processed pressure from \cite{arbogast1995implementation}. This improvement is important in flow coupled with transport problems and it also can be a good candidate for a recovery-based error estimate evaluation. In a recent work \cite{amanbek2018new}, a \textit{posteriori} error analysis was shown for the incompressible flow problems without recovery of velocity.

The remaining part of the paper proceeds as follows. Section \ref{modelformulation} of this paper will describe model formulation with different view of EVMFEM. In Section \ref{method}, the proposed numerical method will be discussed. Section \ref{numericalexample} shows numerical results. Section \ref{conclusion} summarizes the results of this work and draws conclusions.  
%------------------------------------------------------------------------------
\section{Model formulation} \label{modelformulation}
We start by giving the model formulation for the incompressible single-phase flow. For the convenience of reader we repeat the relevant material of domain decomposition method, discrete formulation with Enhanced Velocity from \cite{wheeler2002enhanced}. We next describe the proposed different view of Enhanced Velocity Discrete Scheme with projection operator.
\subsection{Governing equations of the incompressible flow}
We consider the incompressible single-phase flow model for pressure $p$ and the Darcy velocity $\textbf{u}$: 
\begin{align} 
\textbf{u} &=-\mathbf{K} \nabla p \qquad \text{in} \quad \Omega,  \label{eq:a} \\ 
\nabla \cdot \textbf{u} &=f \qquad \qquad \text{in} \quad \Omega, \label{eq:b}  \\ 
p &=g \qquad \qquad \text{on} \quad \partial \Omega   \label{eq:c} 
\end{align}
where $\Omega \in \mathbb{R}^d( d=2$ or $3$) is multiblock domain, $f \in L^2(\Omega)$  and  $\mathbf{K}$ is a symmetric, uniformly positive definite tensor representing the permeability divided by the viscosity with $L^{\infty}(\Omega)$ components, for some $0<k_{min}<k_{max} < \infty$
$k_{min} \xi ^T \xi \le \xi^T \mathbf{K}(x) \xi \le k_{max} \xi^T \xi \qquad \forall x \in \Omega \quad \forall \xi \in \mathbb{R}^d$, under the Dirichlet boundary condition.

A weak variational form of the fluid flow problem $(\ref{eq:a})-(\ref{eq:c})$ is to find a pair $\textbf{u} \in \mathbf{V}$, $p \in W$
\begin{align} 
\left(\mathbf{K}^{-1}\textbf{u}, \mathbf{v} \right) -\left(p, \nabla \cdot \mathbf{v} \right) &= - \langle g, \mathbf{v} \cdot {\boldsymbol \nu} \rangle_{\partial \Omega} \qquad & \forall \mathbf{v} \in \textbf{V} \label{eq:2_4} \\ 
\left(\nabla \cdot \textbf{u}, w \right) &=\left(f,w\right) \qquad \qquad &\forall w \in W \label{eq:2_5}  
\end{align}
where ${\boldsymbol \nu}$ is the outward unit normal to $\partial \Omega$, $\textbf{V}$ is $H({\rm div}; \Omega) =\{\mathbf{v}  \in \left(L^2(\Omega)\right)^d: \nabla \cdot \mathbf{v}  \in L^2(\Omega)\} $ and equipped with the norm $\norm{\mathbf{v} }_V=\left(\norm{\mathbf{v} }^2+\norm{\nabla \cdot \mathbf{v} }^2\right)^{\frac{1}{2}}$ and the pressure   the space is 
$W =L^2(\Omega)$
and the corresponding norm  $\norm{w}_W=\norm{w}.$.

%------------------------------------------------------------------------------------------------------------------------------------
\subsubsection*{Discrete formulation}
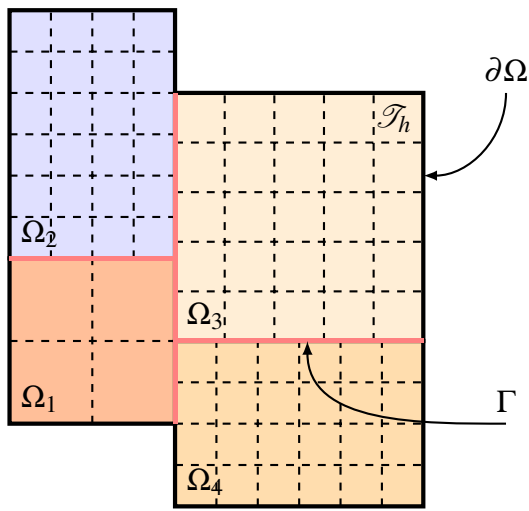
\begin{figure}[H]
	\centering
	\begin{tikzpicture}[thick,scale=1.1] %1.35
	
	\def\xa{-2}
	\def\xb{0}
	\def\xc{3}
	
	\def\ya{-2}
	\def\yb{-1}
	\def\yc{0}
	\def\yd{1}
	\def\ye{3}
	\def\yf{4}
	
	%draw right cube
	\coordinate (A1) at (\xa, \yd);
	\coordinate (A2) at (\xb, \yd);
	\coordinate (A3) at (\xb, \yf);
	\coordinate (A4) at (\xa, \yf);
	\coordinate (B1) at (\xa, \yb);
	\coordinate (B2) at (\xb, \yb);
	%	\coordinate (B3) at (1.7, 1.5);
	%	\coordinate (B4) at (1.7, 0.5);
	\coordinate (C1) at (\xb, \ya );
	\coordinate (C2) at (\xc, \ya);
	\coordinate (C3) at (\xc, \ye);
	\coordinate (C4) at (\xb, \ye);
	%	\coordinate (C5) at (4.5, 2.5);
	%	\coordinate (C6) at (4.5, 4);
	\coordinate (D1) at (\xb, \yc);
	\coordinate (D2) at (\xc, \yc);

	\draw[fill=blue!30,opacity=0.4] (A1) -- (A2) -- (A3) -- (A4)--cycle;
	\draw[fill={rgb:orange,1;yellow,2;pink,5},opacity=0.6] (C1) -- (C2) -- (D2) -- (D1)--cycle;
	\draw[fill={rgb:orange,1;yellow,2;pink,5},opacity=0.3] (D1) -- (D2) -- (C3) -- (C4)--cycle;
	\draw[fill={rgb:orange,1;yellow,1;red,2},opacity=0.4] (B1) -- (B2) -- (A2) -- (A1)--cycle;
	%	\draw[fill={rgb:red,1;green,2;blue,3},opacity=0.4] (D2) -- (D1) -- (B2) -- (A4)--cycle;
	
	% draw grids
	
	\draw[step=5mm,black,thick,dashed] (A1) grid (A3); 
	\draw[step=10mm,black, thick,dashed] (B1) grid (A2); 
	\draw[step=5mm,black,thick,dashed] (C1) grid (D2); 
	\draw[step=6mm,black, thick,dashed] (D1) grid (C3);

	\draw[black,line width=1.5mm,ultra thick] (B1) rectangle (A2);
	\draw[black,line width=1.5mm,ultra thick] (C1) rectangle (D2);
	\draw[black,line width=1.5mm,ultra thick] (D1) rectangle (C3);
	\draw[black,line width=1.5mm,ultra thick] (A1) rectangle (A3);
	%	%nested grid
	%	        \draw[step=5mm, black,ultra thick,dashed] (B1) grid (C2);
	%	        \draw[black,ultra thick] (C4) rectangle (C6);
	
	%filling
	\fill[red!50,ultra thick] (\xb-0.025,\yb) rectangle (\xb+0.025,\ye);%interface
	\fill[red!50,ultra thick] (\xa,\yd-0.025) rectangle (\xb,\yd+0.025);%interface
	\fill[red!50,ultra thick] (\xb,\yc-0.025) rectangle (\xc,\yc+0.025);%interface

	%draw nodes
	\node [above right] at (B1) { $\Omega_1$};
	\node [above right] at (A1) {$\Omega_2$};
	\node [above right] at (D1) {$\Omega_3$};
	\node [above right] at (C1) {$\Omega_4$};
	\node [below left] at (C3) {$\mathcal{T}_h$};

	%labelling
	\draw[-latex,thick](\xc+1,\ye)node[above]{$\partial \Omega$} to[out=270,in=0] (\xc,\ye-1);
	\draw[-latex,thick](\xc+1,\yb)node[above]{$\Gamma$}  to[out=180,in=270] (\xc-1.4,\yc);	
	\end{tikzpicture}
	\caption{Illustration of a domain $\Omega$ with subdomains $\Omega_i$ and non-matching mesh discretization $\mathcal{T}_h$.}
	\label{fig:ddm}
\end{figure}

Let $\Omega$ be decomposed into non-overlapping small subdomains, see Fig. \ref{fig:ddm}. We consider 
\begin{align*}
\Omega = \left( \bigcup_{i=1}^{N_b} \bar{\Omega}_i \right)^{o}, \; \Gamma_{i,j}=\partial \Omega_i \bigcap \partial \Omega_j, \; \Gamma = \left( \bigcup^{N_b}_{i,j=1} \bar{\Gamma}_{i,j} \right)^{o}, \; \Gamma_i = \Omega_i \bigcap \Gamma=\partial \Omega_j \setminus \partial \Omega.
\end{align*}
This implies that the domain is divided into $N_b$ subdomains, the interface between $i^{th}$ and $j^{th}$ subdomains($i \ne j$), the interior subdomain interface for $i^{th}$ subdomain and union of all such interfaces, respectively. 
Let $\mathcal{T}_{h,i}$ be a conforming, quasi-uniform and rectangular partition of $\Omega_i$, $1 \le i \le N_b$, with maximal element diameter $h_i$. We then set $\mathcal{T}_{h}=\cup_{i=1}^{n} \mathcal{T}_{h,i} $ and denote $h$ the maximal element diameter in $\mathcal{T}_h$; note that $\mathcal{T}_h$ can be nonmatching as neighboring meshes $\mathcal{T}_{h,i}$ and $\mathcal{T}_{h,j}$ need not match on $\Gamma_{i,j}$. We assume that all mesh families are shape-regular.

In Enhanced Velocity scheme setting, the velocity basis functions are based on the traditional Raviart-Thomas spaces of lowest order on rectangles for $d=2$ and bricks for $d=3$. The $RT_0$ spaces are defined for any element $T \in \mathcal{T}_{h}$ by the following spaces:
\begin{align*}
&\mathbf{V}_h(T)=\{\mathbf{v}=(v_1,v_2) \; \text{or} \; \mathbf{v}=(v_1,v_2,v_3): v_l=\alpha_l+\beta_l x_l: \alpha_l, \beta_l \in \mathbb{R};l=1,..d \}, \\
&W_h(T) =\{w=\text{constant}\}.
\end{align*}
The pressure finite element approximation space on $\Omega$ is taken to be as
$ W_h(\Omega) =\{w\in L^2(\Omega): w \biggr\rvert_E \in W_h(T), \forall T \in \mathcal{T}_h \}.$
In addition, a vector function in $\mathbf{V}_h$ can be determined uniquely by its normal components $\mathbf{v} \cdot \nu $ at midpoints of edges (in 2D) or face (in 3D) of $T$. 
%This is illustrated in Figure \ref{fig:RT0}.
The degrees of freedom of $\mathbf{v} \in \mathbf{V}_h(T)$ were created by these normal components. The degree of freedom for a pressure function $p \in W_h(T)$ is at center of $T$ and piecewise constant inside of $T$. 

Let us formulate $RT_0$ space on each subdomain $\Omega_i$ for partition $\mathcal{T}_h$
\begin{align*}
\mathbf{V}_{h,i}=\{\mathbf{v} \in H({\rm div};\Omega_i): \mathbf{v} \biggr\rvert_T \in \mathbf{V}_h(T), \forall T \in \mathcal{T}_{h,i} \}  \qquad i \in \{1,...n\}
\end{align*}
and then
\begin{align*}
\mathbf{V}_{h}= \bigoplus_{i=1}^n \mathbf{V}_{h,i}.
\end{align*}
Although the normal components of vectors in $\mathbf{V}_h$ are continuous between elements within each subdomains, the reader may see $\mathbf{V}_h$ is \textit{not} a subspace of $H({\rm div}; \Omega)$, because the normal components of the velocity vector may not match on subdomain interface $\Gamma$. 

\begin{figure}[!tbp]
	\centering	
	\begin{tikzpicture}[thick,scale=1.0, dot/.style = {outer sep = +0pt, inner sep = +0pt, shape = circle, draw = black, label = {#1}},
	small dot/.style = {minimum size = 1pt, dot = {#1}},
	big dot/.style = {minimum size = 8pt, dot = {#1}},
	line join = round, line cap = round, >=triangle 45
	]
	
	%draw element rectangle
	\coordinate (A1) at (0, -1);
	\coordinate (A2) at (0, 1);
	\coordinate (A3) at (2, 1);
	\coordinate (A4) at (2, -1);
	
	\coordinate (B1) at (-2.6,-2);
	\coordinate (B2) at (-2.6, 1.5);
	\coordinate (B3) at (0,1.5);
	\coordinate (B4) at (0,-2);
	
	\coordinate (C1) at (-2.6,0);
	\coordinate (C2) at (0, 0);
	\coordinate (C3) at (2,0);
	
	\coordinate (D1) at (-2.6, -1);
	\coordinate (D2) at (-2.6,1);

	%filling
	\fill[fill={rgb:orange,1;yellow,2;pink,5}, opacity=0.2] (A1) rectangle (A3); %fill={rgb:orange,1;yellow,2;pink,5},opacity=0.3
	
	%rectangles
	\draw[very thick] (A1) -- (A2) -- (A3) -- (A4) -- cycle;
	
	\draw[fill={rgb:orange,1;yellow,2;blue,2},opacity=0.2, very thick] (B1) -- (B4)-- (C2)--(C1);
	\draw[fill={rgb:orange,2;yellow,1;green,1},opacity=0.2, very thick] (C1) -- (C2)-- (B3)--(B2);
	
	\draw[very thick] (B4) -- (B3);
	\draw[very thick] (B3) -- (B2);
	\draw[very thick] (C1) -- (C2);
	\draw[dotted, color =red] (C2) -- (C3);
	\draw[dotted, color =red] (A2) -- (D2);
	\draw[dotted, color =red] (A1) -- (D1);
	\draw[very thick] (B1) -- (B4);
	\draw[very thick] (B1) -- (B2);

	%velocity dof
	\node[font = \large, color =red] (e1)  at (0, 0.5) { \textbf{$\times$}};
	\node[font = \large, color =red] (e2)  at (0,-0.5) {$\times$};	
	\node[font = \large ] (e3)  at (1, -1) {$\times$};
	\node[font = \large] (e4)  at (2, 0) {$\times$};
	\node[font = \large] (e5)  at (1, 1) {$\times$};
	
	%labeling
	\node [above, font = \large] at (0,-2.7) { $\Gamma_{i,j}$};
	\node [above left, font = \large] at (0, 0.5) { $e_1$};
	\node [above left, font = \large] at (0,-0.5) { $e_2$};
	\node [right, font = \large] at (1, 0.5) { $T_1$};
	\node [right, font = \large] at (1,-0.5) { $T_2$};
	
	\end{tikzpicture}
	\captionof{figure}{Degrees of freedom for the Enhanced Velocity space.}
	\label{fig:EVMFEM_Gamma}
\end{figure}
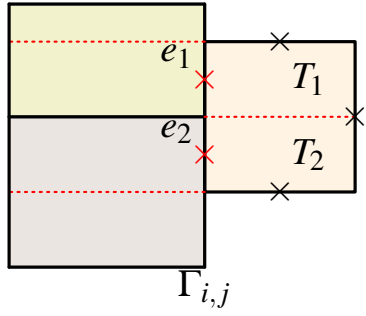

Let us define $\mathcal{T}_{h,i,j}$ as the intersection of the traces of $\mathcal{T}_{h,i}$ and $\mathcal{T}_{h,j}$, and let $\mathcal{T}^{\Gamma}_h=\bigcup_{1 \le i \le j \le N_b} \mathcal{T}_{h,i,j}$. We require that $\mathcal{T}_{h,i}$ and $\mathcal{T}_{h,j}$ need to align with the coordinate axes. Fluxes are constructed to match on each element $e \in \mathcal{T}^{\Gamma}_h$. We consider any element $T \in \mathcal{T}_{h,i}$ that shares at least one edge with the interface $\Gamma$, i.e., $T \cap \Gamma_{i,j} \ne \emptyset$, where $1 \le i,j \le N_b$ and $i \ne j$. Then newly defined interface grid introduces a partition of the edge of $T$. This partition may be extended into the element $T$ as shown in Fig. \ref{fig:EVMFEM_Gamma}.

Such partitioning helps to construct fine-scale velocities that is in $H(\textbf{ div}, \Omega)$. So we represent a basis function $\mathbf{v}_{T_k}$ in the $\mathbf{V}_h(T_k)$ space ($RT_0$) for given $T_k$ with the following way: 
\begin{align*}
\mathbf{v}_{T_k} \cdot \nu =
\begin{cases*}
1, \qquad {\rm on} \; e_k \\
0, \qquad \rm other \; edges
\end{cases*}
\end{align*}
i.e. a normal component $\mathbf{v}_{T_k} \cdot \nu $ equal to one on $e_k$ and zero on all other edges(faces) of $T_k$.
Let $\mathbf{V}^{\Gamma}_h$ be span of all such basis functions defined on all sub-elements induced the interface discretization $\mathcal{T}_{h, i,j}$. Thus, the enhanced velocity space $\mathbf{V}^*_h$ is taken to be as
\begin{align*}
\mathbf{V}^*_{h}= \bigoplus_{i=1}^n \mathbf{V}^0_{h,i} \bigoplus \mathbf{V}^{\Gamma}_{h} \cap H({\rm div}; \Omega).
\end{align*}
where $\mathbf{V}^0_{h,i} = \{ \mathbf{v} \in \mathbf{V}_{h,i} : \mathbf{v} \cdot \nu =0  \text{  on  } \Gamma_{i}\}$ is the subspace of $\mathbf{V}_{h,i}$. The finer grid velocity allows to velocity approximation on the interface and then form the $H({\rm div}, \Omega)$ conforming velocity space. Some difficulties arise, however, in analysis of method and implementation of robust linear solver for such modification of $RT_0$ velocity space at all elements, which are adjacent to the interface $\Gamma$.
We now formulate the discrete variational form of equations $(\ref{eq:a})-(\ref{eq:c})$ as: Find  $\mathbf{u}_h \in \mathbf{V}^*_{h} $ and $p_h \in W_h$ such that
\begin{align} 
\left(K^{-1}\textbf{u}_h, \mathbf{v} \right) &=\left(p_h, \nabla \cdot \mathbf{v} \right)- \langle g, \mathbf{v} \cdot {\boldsymbol \nu} \rangle_{\partial \Omega} \qquad & \forall \mathbf{v} \in \textbf{V}^*_h  \label{eq:EV_eq1} \\ 
\left(\nabla \cdot \textbf{u}_h, w \right) &=\left(f,w\right) \qquad \qquad &\forall w \in W_h\label{eq:EV_eq2}  
\end{align}

%------------------------------------------------------------------------------------------------------------------------------------

\subsection{A different view of the EVMFEM in the Discrete Variational Formulation} \label{sec:diffEV}

We consider the discrete variational form that is given in (\ref{eq:EV_eq1})-(\ref{eq:EV_eq2}).
Find  $\mathbf{u}_h \in \mathbf{V}^*_{h} $ and $p_h \in W_h$ such that
\begin{align} 
\left(K^{-1}\textbf{u}_h, \mathbf{v} \right)_{M, T} &=\left(p_h, \nabla \cdot \mathbf{v} \right)- \langle g, \mathbf{v} \cdot {\boldsymbol \nu} \rangle_{\partial \Omega} \qquad & \forall \mathbf{v} \in \textbf{V}^*_h  \label{eq:EV_h1} \\ 
\left(\nabla \cdot \textbf{u}_h, w \right) &=\left(f,w\right) \qquad \qquad &\forall w \in W_h\label{eq:EV_h2}  
\end{align}

We exploit the approximation inner product and for $\mathbf{v}, \mathbf{q} \in \mathbb{R}^d$
\begin{align*}
\left( \mathbf{v}, \mathbf{q} \right)_{M,T} = 
\begin{cases*}
\left( v_{x}, q_y \right)_{T_x, M_y} + \left( v_{y}, q_y \right)_{M_x,T_y} \quad &\text{if } d =2, \\
\left( v_{x}, q_y \right)_{T_x, M_y, M_z} + \left( v_{y}, q_y \right)_{M_x, T_y, M_z}+\left( v_{z}, q_z \right)_{M_x,M_y,T_z} \quad &\text{if } d=3.
\end{cases*}
\end{align*}
where $T_{(\cdot)}$ and $M_{(\cdot)}$ denote the the trapezoidal  and midpoint quadrature rules in each coordinate direction respectively, see \cite{russell1983finite}. In particularly, we take $\mathbf{v}=\mathbf{K}^{-1}\mathbf{u}_h$ and $\mathbf{q}=\mathbf{v}$.

It is easily proven that the finite variational form (\ref{eq:EV_h1})-(\ref{eq:EV_h2}) is equivalent to finding $\mathbf{u}_h \in \mathbf{V}^*_h$, $p_h \in W_h$, $1 \le i \le N_b$, such that
\begin{align} 
\left(\mathbf{K}^{-1}\mathbf{u}_h, \mathbf{v} \right)_{\Omega_i, M, T} - \left(p_h, \nabla \cdot \mathbf{v} \right)_{\Omega_i} &=- \langle g, \mathbf{v} \cdot {\boldsymbol \nu} \rangle_{\partial \Omega_i \cap \Gamma_D} \quad & \forall \mathbf{v} \in \mathbf{V}^0_{h,i} \label{eq:6_8p} \\ 
\left(\nabla \cdot \mathbf{u}_h, w \right)_{\Omega_i} &=\left(f,w\right)_{\Omega_i} \quad \qquad &\forall w \in W_{h,i} \label{eq:6_9p}  \\
\sum_{i=1}^{N_b} \{\left(\mathbf{K}^{-1}\mathbf{u}_h, \mathbf{v}^{EV}\right)_{\Omega_i, M, T}-\left(p_h, \nabla \cdot \mathbf{v}^{EV}\right)_{\Omega_i}\} &= 0 \quad & \forall \mathbf{v}^{EV} \in \textbf{V}^{\Gamma}
\end{align}
We note that similar the discrete variational formulation was proposed in \cite{glowinski1988domain} with conjugate gradient method.
We want to share the idea for small number of discretization elements that can be applied for a large number of elements. Thus, we consider two subdomains,i.e., $\Omega =\bar{\Omega}_1 \cup \bar{\Omega}_2 $ and $\Gamma$ is the interface. Then
\begin{equation*}
\mathbf{V}^*_{h}=\left( \mathbf{V}^0_{h,1} \oplus \mathbf{V}^0_{h,2} \oplus \mathbf{V}^{\Gamma}_{h} \right) %\cap H({\rm div}; \Omega) 
\end{equation*}
Consider equations 
\begin{align} 
\left(\mathbf{K}^{-1}\textbf{u}_h, \mathbf{v} \right)_{M, T} &=\left(p_h, \nabla \cdot \mathbf{v} \right) \qquad & \forall \mathbf{v} \in \mathbf{V}^{\Gamma}_{h} \label{eq:4_9}
\end{align}

\begin{figure}[!tbp]
	\centering	
	\begin{tikzpicture}[thick,scale=0.7, dot/.style = {outer sep = +0pt, inner sep = +0pt, shape = circle, draw = black, label = {#1}},
	small dot/.style = {minimum size = 1pt, dot = {#1}},
	big dot/.style = {minimum size = 8pt, dot = {#1}},
	line join = round, line cap = round, >=triangle 45
	]
	
	\def\xa{0}
	\def\xb{6}
	\def\xc{12}
	
	\def\ya{0}
	\def\yb{6}
	
	\def\h{2}
	\def\hr{3}
	
	%draw right cube
	\coordinate (A1) at (\xa, \ya);
	\coordinate (A2) at (\xb, \ya);
	\coordinate (A3) at (\xb, \yb);
	\coordinate (A4) at (\xa, \yb);
	\coordinate (B1) at (\xc, \ya);
	\coordinate (B2) at (\xc, \yb);
	
	%draw rectangle
	\draw[] (A1) -- (A2) -- (A3) -- (A4)--cycle;
	\draw[] (A2) -- (B1) -- (B2) -- (A3)--cycle;
	
	% draw grids
	\draw[step=20mm,black, thin, dashed] (A1) grid (A3); 
	\draw[step=30mm,black, thin, dashed] (A2) grid (B2); 
	
	%fill subelements
	\fill[fill={rgb:orange,1;yellow,2;blue,2},opacity=0.2, very thick] (\xb-\h,\ya) rectangle (\xb,\ya+\h);
	\fill[fill={rgb:orange,2;yellow,1;green,1},opacity=0.2, very thick] (\xb-\h,\ya+\h) rectangle (\xb,\ya+2*\h);
	\fill[fill={rgb:orange,4;yellow,1;green,0},opacity=0.4, very thick] (\xb-\h,\ya+2*\h) rectangle (\xb,\ya+3*\h);
	\fill[fill=brown!40,opacity=0.6, very thick] (\xb,\ya) rectangle (\xb+\hr,\ya+\hr);
	\fill[fill={rgb:orange,1;yellow,2;pink,5}, opacity=0.2, very thick] (\xb,\ya+\hr) rectangle (\xb+\hr,\ya+2*\hr);
	
	%velocity dof
	\node[font = \Large, color =red] (e1)  at (\xb, 1) { $\times$};
	\node[font = \Large, color =red] (e2)  at (\xb, 2.5) {$\times$};
	\node[font = \Large, color =red] (e3)  at (\xb, 3.5) { $\times$};
	\node[font = \Large, color =red] (e4)  at (\xb, 5) {$\times$};	
	
	%pressure dof
	\node[ fill = green!40, big dot] (p)  at (\xb-1, \ya+1) {};
	\node[ fill = green!40, big dot] (p)  at (\xb-1, \ya+3) {};
	\node[ fill = green!40, big dot] (p)  at (\xb-1, \ya+5) {};
	\node[ fill = green!40, big dot] (p)  at (\xb+1.5, \ya+1.5) {};
	\node[ fill = green!40, big dot] (p)  at (\xb+1.5, \ya+4.5) {};

	%labeling
	\node [below, font = \large] at (\xb,\ya) { $\Gamma_{i,j}$};
	\node [above right, font = \large] at (\xa,\ya) { $\Omega_{L}$};
	\node [above left, font = \large] at (\xc,\ya) { $\Omega_{R}$};
	\node [below, font = \large,color=black!50!green] at (\xb-.5*\h,\ya) { $p_{L}$};
	\node [below, font = \large, color=black!50!green] at (\xb+.5*\hr,\ya) { $p_{R}$};
	\end{tikzpicture}
	\caption{The spatial domain and illustration of Enhanced Velocity values on the interface}
	\label{fig:EVMFEM_two_domains_flux}
\end{figure}
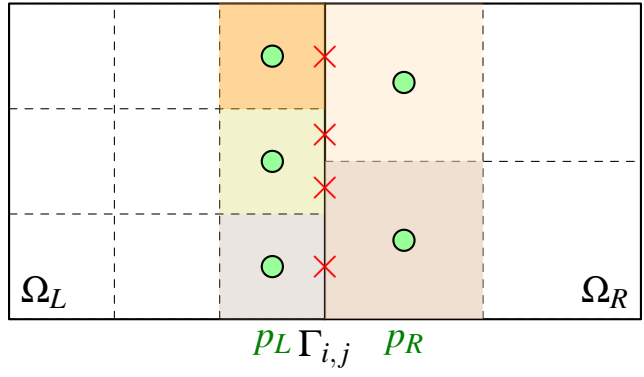
These allow us to express $\mathbf{u}^{\Gamma}_h$ in terms of the one-element layers along $\Gamma$, it is shown in Fig. \ref{fig:EVMFEM_two_domains_flux}:
\begin{equation}
\mathbf{u}^{\Gamma}_h=A_1p_L+A_2p_R \label{eq:projection_pressure}
\end{equation}
Now we consider each subdomain separately with ghost layers. We define $L_2$-projection of Enhanced Velocity space at interface $\Gamma_{i,j}$ to each subdomain space $\partial \Omega_i \cap \Gamma_{i,j}$ such that $\mathcal{P}^{i}_h : \mathbf{V}^*_h \rightarrow \mathbf{V}_{h,i}$. 
%denote $\mathcal{P}^{i}_h$ be the $\Lambda_H(\Gamma)$-orthogonal projection onto $\Lambda_H$,
\begin{align*}
\mathcal{P}^{i}_h :\mathbf{V}^*_h(\Gamma) \to \mathbf{V}_{h,i}(\Gamma_i) \qquad  for \; \psi \in L^2(\Gamma), \qquad \langle \left( \psi  - \mathcal{P}^{i}_h \psi\right) \cdot {\boldsymbol \nu}_i, \mathbf{v} \cdot {\boldsymbol \nu}_i \rangle_{\Gamma}=0 \qquad \forall \mathbf{v} \in \mathbf{V}_{h,i}.
\end{align*}

We denote 
\begin{align*}
\mathbf{u}^{\Gamma}_{h,i}=\mathcal{P}^{i}_h \mathbf{u}^{\Gamma}_h, \qquad i=L \; or \; R
\end{align*}

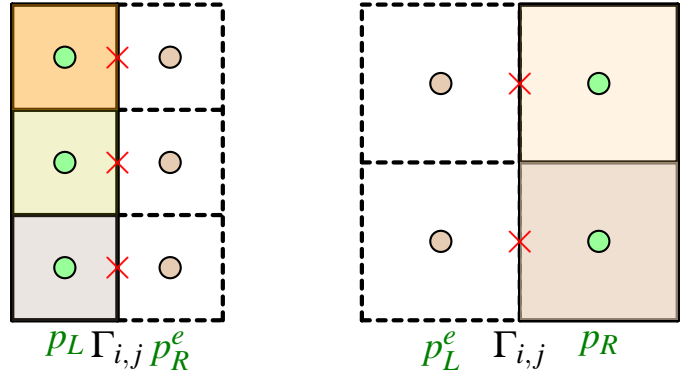
\begin{figure}[!tbp]
	\centering	
	\begin{tikzpicture}[thick,scale=0.7, dot/.style = {outer sep = +0pt, inner sep = +0pt, shape = circle, draw = black, label = {#1}},
	small dot/.style = {minimum size = 1pt, dot = {#1}},
	big dot/.style = {minimum size = 8pt, dot = {#1}},
	line join = round, line cap = round, >=triangle 45
	]
	
	\def\xa{0}
	\def\xb{2}
	\def\xc{4}
	
	\def\ya{0}
	\def\yb{6}
	
	\def\h{2}
	\def\hr{3}
	
	%draw right cube
	\coordinate (A1) at (\xa, \ya);
	\coordinate (A2) at (\xb, \ya);
	\coordinate (A3) at (\xb, \yb);
	\coordinate (A4) at (\xa, \yb);
	\coordinate (B1) at (\xc, \ya);
	\coordinate (B2) at (\xc, \yb);

	% draw grids
	\draw[step=20mm,black,ultra thick] (A1) grid (A3); 
	\draw[step=20mm,black,ultra thick,dashed] (A2) grid (\xb+\h,\yb);
	
	%draw rectangle
	\draw[line width=0.5mm] (A1) -- (A2) -- (A3) -- (A4)--cycle;

	%fill subelements
	\fill[fill={rgb:orange,1;yellow,2;blue,2},opacity=0.2, very thick] (\xb-\h,\ya) rectangle (\xb,\ya+\h);
	\fill[fill={rgb:orange,2;yellow,1;green,1},opacity=0.2, very thick] (\xb-\h,\ya+\h) rectangle (\xb,\ya+2*\h);
	\fill[fill={rgb:orange,4;yellow,1;green,0},opacity=0.4, very thick] (\xb-\h,\ya+2*\h) rectangle (\xb,\ya+3*\h);
	
	%velocity dof
	\node[font = \Large, color =red] (e1)  at (\xb, 1) { $\times$};
	\node[font = \Large, color =red] (e2)  at (\xb, 3) {$\times$};
	\node[font = \Large, color =red] (e4)  at (\xb, 5) {$\times$};	
	
	%pressure dof
	\node[ fill = green!40, big dot] (p)  at (\xb-.5*\h, \ya+.5*\h) {};
	\node[ fill = green!40, big dot] (p)  at (\xb-.5*\h, \ya+1.5*\h) {};
	\node[ fill = green!40, big dot] (p)  at (\xb-.5*\h, \ya+2.5*\h) {};
	\node[ fill = brown!40, big dot] (p)  at (\xb+.5*\h, \ya+.5*\h) {};
	\node[ fill = brown!40, big dot] (p)  at (\xb+.5*\h, \ya+1.5*\h) {};
	\node[ fill = brown!40, big dot] (p)  at (\xb+.5*\h, \ya+2.5*\h) {};

	%labeling
	\node [below, font = \large] at (\xb,\ya) { $\Gamma_{i,j}$};
	\node [below, font = \large,color=black!50!green] at (\xb-.5*\h,\ya) { $p_{L}$};
	\node [below, font = \large, color=black!50!green] at (\xb+.5*\h,\ya) { $p^e_{R}$};
	\end{tikzpicture}
	\qquad \qquad
	\begin{tikzpicture}[thick,scale=0.7, dot/.style = {outer sep = +0pt, inner sep = +0pt, shape = circle, draw = black, label = {#1}},
	small dot/.style = {minimum size = 1pt, dot = {#1}},
	big dot/.style = {minimum size = 8pt, dot = {#1}},
	line join = round, line cap = round, >=triangle 45
	]
	
	\def\xa{0}
	\def\xb{3}
	\def\xc{6}
	
	\def\ya{0}
	\def\yb{6}
	
	\def\h{2}
	\def\hr{3}
	
	%draw right cube
	\coordinate (A1) at (\xa, \ya);
	\coordinate (A2) at (\xb, \ya);
	\coordinate (A3) at (\xb, \yb);
	\coordinate (A4) at (\xa, \yb);
	\coordinate (B1) at (\xc, \ya);
	\coordinate (B2) at (\xc, \yb);

	% draw grids
	%	\draw[step=30mm,black,ultra thick,dashed] (\xb-\hr,\ya) grid (A3); 
	
	% draw grids
	\draw[step=30mm,black,ultra thick] (A2) grid (B2); 
	\draw[step=30mm,black,ultra thick,dashed] (\xb-\hr,\ya) grid (A3); 
	
	%draw rectangle
	\draw[line width=0.5mm] (A2) -- (B1) -- (B2) -- (A3)--cycle;
	
	%fill subelements
	\fill[fill=brown!40,opacity=0.6, very thick] (\xb,\ya) rectangle (\xb+\hr,\ya+\hr);
	\fill[fill={rgb:orange,1;yellow,2;pink,5}, opacity=0.2, very thick] (\xb,\ya+\hr) rectangle (\xb+\hr,\ya+2*\hr);
	
	%velocity dof
	\node[font = \Large, color =red] (e1)  at (\xb, 1.5) { $\times$};
	\node[font = \Large, color =red] (e2)  at (\xb, 4.5) {$\times$};
	
	%pressure dof
	%	\node[ fill = green!40, big dot] (p)  at (\xb-.5*\h, \ya+.5*\h) {};
	%	\node[ fill = green!40, big dot] (p)  at (\xb-.5*\h, \ya+1.5*\h) {};
	%	\node[ fill = green!40, big dot] (p)  at (\xb-.5*\h, \ya+2.5*\h) {};
	%	\node[ fill = brown!40, big dot] (p)  at (\xb+.5*\h, \ya+.5*\h) {};
	\node[ fill = green!40, big dot] (p)  at (\xb+1.5, \ya+1.5) {};
	\node[ fill = green!40, big dot] (p)  at (\xb+1.5, \ya+4.5) {};
	\node[ fill = brown!40, big dot] (p)  at (\xb-.5*\hr, \ya+.5*\hr) {};
	\node[ fill = brown!40, big dot] (p)  at (\xb-.5*\hr, \ya+1.5*\hr) {};

	%labeling
	\node [below, font = \large] at (\xb,\ya) { $\Gamma_{i,j}$};
	\node [below, font = \large,color=black!50!green] at (\xb-.5*\hr,\ya) { $p^e_{L}$};
	\node [below, font = \large, color=black!50!green] at (\xb+.5*\hr,\ya) { $p_{R}$};
	\end{tikzpicture}
	
	\caption{Example of left ($\Omega_L$) and right ($\Omega_R$) domains with ghost layers ghost layers}
	\label{fig:domainghostlayer}
\end{figure}

In subdomain $\Omega_i$, we define $p^e_i$ in the following way
\begin{align} 
\left(\mathbf{K}^{-1}\widetilde{\mathbf{u}}_h, \mathbf{v} \right)_{M, T, \Omega_i} &=\left(p_h, \nabla \cdot \mathbf{v} \right)_{\Omega_i}- \langle p^{e}_{i}, \mathbf{v} \cdot {\boldsymbol \nu} \rangle_{\Gamma} \qquad & \forall \mathbf{v} \in \mathbf{V}^{\Gamma}_{h, i} \text{  s.t.} \mathbf{v}\cdot {\boldsymbol \nu} =0 \text{ on } \partial \Omega^*_i \label{eq:5_6}
\end{align}
where $\Omega^*$ is union of all elements $T$ that shares edge (2D) or face (3D) with $\Gamma_i$ and $p^e_i$ ghost layers pressure values, and $\widetilde{\mathbf{u}}_h=\mathcal{P}^{i}_h(\mathbf{u}_h)$. Such ghost layers are depicted in the Fig. \ref{fig:domainghostlayer}.  
Then, for $i=L$, we have 
\begin{equation}
\textbf{u}^{\Gamma}_{h,L}=A^L_1 p_L +A^L_2 p^e_L  \label{eq:projection_subdomain}
\end{equation}

We compare equation (\ref{eq:projection_subdomain}) and the projected to $\Omega_i$ pressure equation (\ref{eq:projection_pressure}):
\begin{equation}
\mathcal{P}^{L}_h \textbf{u}^{\Gamma}_h=\mathcal{P}^{L}_h A_1 p_L +\mathcal{P}^{L}_h A_2 p_R \\
\end{equation}
Since $\mathbf{u}^{\Gamma}_{h,L}=\mathcal{P}^{L}_h \textbf{u}^{\Gamma}_h$, $A^L_1=\mathcal{P}^{L}_h A_1$ , we have the following
\begin{equation}
A^L_2 p^e_L = \mathcal{P}^{L}_h  A_2 p_R
\end{equation}
$A^L_2$ is non-singular and diagonal matrix, since $\mathbf{K}$ is SPD.
\begin{equation}
p^e_L = \left(A^L_2\right)^{-1} \mathcal{P}^{L}_h A_2 p_R.
\end{equation}
Similarly, we can obtain 
\begin{equation}
p^e_R = \left(A^R_1\right)^{-1} \mathcal{P}^{R}_h  A_1 p_L
\end{equation}

In non-linear problems including slightly compressible flow or multiphase flow in heterogeneous porous media, this approach could be applied analogously by taking into account ghost layers values arising from $e_i$ in each Newton iteration. So during Block Jacobi iteration variables $p^{e,k-1}_{L}$, $p^{e, k-1}_{R}$ is computed by utilizing given $p^{k-1}_{L}$, $p^{k-1}_{R}$ and then solve decoupled subdomain problems with Dirichlet boundary conditions  $p^{e, k-1}_{i}$, $i=L,R$ to find $u^{k},p^{k}$.
%------------------------------------------------------------------------------
\section{Methods} \label{method}
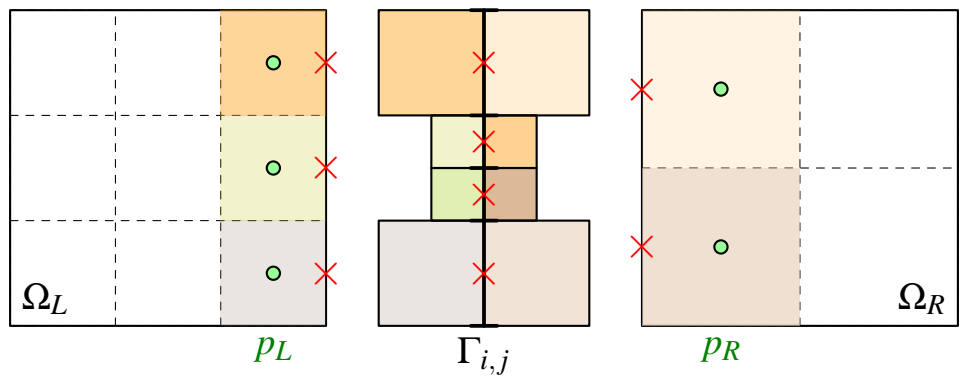
\begin{figure}[H]
	\centering	
	\begin{tikzpicture}[thick,scale=0.7, dot/.style = {outer sep = +0pt, inner sep = +0pt, shape = circle, draw = black, label = {#1}},
	small dot/.style = {minimum size = 1pt, dot = {#1}},
	big dot/.style = {minimum size = 5pt, dot = {#1}},
	line join = round, line cap = round, >=triangle 45
	]
	
	\def\xa{0}
	\def\xb{6}
	\def\xi{9}
	\def\xc{12}
	\def\xd{18}
	
	\def\ya{0}
	\def\yb{6}
	
	\def\h{2}
	\def\hr{3}
	
	%draw right cube
	\coordinate (A1) at (\xa, \ya);
	\coordinate (A2) at (\xb, \ya);
	\coordinate (A3) at (\xb, \yb);
	\coordinate (A4) at (\xa, \yb);
	\coordinate (B1) at (\xd, \ya);
	\coordinate (B2) at (\xd, \yb);
	\coordinate (B3) at (\xc, \ya);
	\coordinate (B4) at (\xc, \yb);
	\coordinate (C1) at (\xi-\h, \yb-\h);
	\coordinate (C2) at (\xi+\h, \yb);
	\coordinate (D1) at (\xi-0.5*\h, \yb-1.5*\h);
	\coordinate (D2) at (\xi+0.5*\h, \yb-\h);
	\coordinate (E1) at (\xi-0.5*\h, \yb-2*\h);
	\coordinate (E2) at (\xi+0.5*\h, \yb-1.5*\h);
	\coordinate (F1) at (\xi-\h, \ya);
	\coordinate (F2) at (\xi+\h, \ya+\h);
	
	%draw rectangle
	\draw[] (A1) -- (A2) -- (A3) -- (A4)--cycle;
	\draw[] (B1) -- (B2) -- (B4) -- (B3) -- cycle;

	% draw grids
	\draw[step=20mm,black, thin,dashed] (A1) grid (A3); 
	\draw[step=30mm,black, thin,dashed] (B3) grid (B2); 
	
	%fill subelements
	\fill[fill={rgb:orange,1;yellow,2;blue,2},opacity=0.2, very thick] (\xb-\h,\ya) rectangle (\xb,\ya+\h);
	\fill[fill={rgb:orange,2;yellow,1;green,1},opacity=0.2, very thick] (\xb-\h,\ya+\h) rectangle (\xb,\ya+2*\h);
	\fill[fill={rgb:orange,4;yellow,1;green,0},opacity=0.4, very thick] (\xb-\h,\ya+2*\h) rectangle (\xb,\ya+3*\h);
	\fill[fill=brown!40,opacity=0.6, very thick] (\xc,\ya) rectangle (\xc+\hr,\ya+\hr);
	\fill[fill={rgb:orange,1;yellow,2;pink,5}, opacity=0.2, very thick] (\xc,\ya+\hr) rectangle (\xc+\hr,\ya+2*\hr);
	\fill[fill={rgb:orange,4;yellow,1;green,0},opacity=0.4, very thick] (C1) rectangle (\xi,\yb);
	\fill[fill={rgb:orange,1;yellow,2;blue,2},opacity=0.2, very thick] (F1) rectangle (\xi,\ya+\h);
	\fill[fill={rgb:orange,2;yellow,1;green,2},opacity=0.3, very thick] (E1) rectangle (\xi,\ya+1.5*\h);
	\fill[fill={rgb:orange,2;yellow,1;green,1},opacity=0.2, very thick] (D1) rectangle (\xi,\ya+2*\h);
	\fill[fill=brown!40,opacity=0.55, very thick] (\xi,\ya) rectangle (F2);
	\fill[fill=brown!60,opacity=0.9, very thick] (\xi,\ya+\h) rectangle (E2);
	\fill[fill={rgb:orange,1;yellow,2;pink,5}, opacity=0.8, very thick] (\xi,\ya+\hr) rectangle (D2);
	\fill[fill={rgb:orange,1;yellow,2;pink,5}, opacity=0.3, very thick] (\xi,\yb-\h) rectangle (C2);

	\draw[draw=black] (C1) rectangle (C2);
	\draw[draw=black] (D1) rectangle (D2);
	\draw[draw=black] (E1) rectangle (E2);
	\draw[draw=black] (F1) rectangle (F2);
	%Interface
	\draw[-, line width=0.5mm] (\xi-0.25,\ya) -- (\xi+0.25,\ya);
	\draw[-, line width=0.5mm] (\xi-0.25,\yb) -- (\xi+0.25,\yb);
	\draw[-, line width=0.5mm] (\xi-0.25,2.0) -- (\xi+0.25,2.0);
	\draw[-, line width=0.5mm] (\xi-0.25,3.0) -- (\xi+0.25,3.0);
	\draw[-, line width=0.5mm] (\xi-0.25,4.0) -- (\xi+0.25,4.0);
	\draw[-, line width=0.5mm] (\xi,\ya) -- (\xi,\yb);
	
	%velocity dof
	\node[font = \Large, color =red] (e1)  at (\xi, 1) { $\times$};
	\node[font = \Large, color =red] (e2)  at (\xi, 2.5) {$\times$};
	\node[font = \Large, color =red] (e3)  at (\xi, 3.5) { $\times$};
	\node[font = \Large, color =red] (e4)  at (\xi, 5) {$\times$};
	\node[font = \Large, color =red] (e1)  at (\xb, 1) { $\times$};
	\node[font = \Large, color =red] (e2)  at (\xb, 3) {$\times$};
	\node[font = \Large, color =red] (e4)  at (\xb, 5) {$\times$};	
	\node[font = \Large, color =red] (e2)  at (\xc, 1.5) {$\times$};
	\node[font = \Large, color =red] (e4)  at (\xc, 4.5) {$\times$};	
	
	%pressure dof
	\node[ fill = green!40, big dot] (p)  at (\xb-1, \ya+1) {};
	\node[ fill = green!40, big dot] (p)  at (\xb-1, \ya+3) {};
	\node[ fill = green!40, big dot] (p)  at (\xb-1, \ya+5) {};
	\node[ fill = green!40, big dot] (p)  at (\xc+1.5, \ya+1.5) {};
	\node[ fill = green!40, big dot] (p)  at (\xc+1.5, \ya+4.5) {};
	
	%labeling
	\node [below, font = \large] at (\xi,\ya) { $\Gamma_{i,j}$};
	\node [above right, font = \large] at (\xa,\ya) { $\Omega_{L}$};
	\node [above left, font = \large] at (\xd,\ya) { $\Omega_{R}$};
	\node [below, font = \large,color=black!50!green] at (\xb-.5*\h,\ya) { $p_{L}$};
	\node [below, font = \large, color=black!50!green] at (\xc+.5*\hr,\ya) { $p_{R}$};
	\end{tikzpicture}
	\caption{The illustration of the velocity improvement at the interface using postprocessing.}
	\label{fig:EVMFEM_postprocess_domains_projection}
\end{figure}

We use the postprocessing procedure associated to pressure and velocity. We first apply locally postprocessing algorithm for given pressure $p_h$ and velocity $\mathbf{u}_h$ which was previously proposed in \cite{arbogast1995implementation} and then Oswald interpolation operator \cite{pencheva2013robust,ainsworth2005robust,karakashian2003posteriori,vohralik2010unified} to have better pressure values. At the interface, we use two-point flux computation method in order to have better approximation of pressure. As a result, the Enhanced Velocity scheme solution of flux can be improved by using a post-processed pressure. 
The key idea is illustrated in Fig. \ref{fig:EVMFEM_postprocess_domains_projection} for resulting approximation of EV scheme that is shown in  Fig. \ref{fig:EVMFEM_two_domains_flux}. 

The velocity at the edge or face is computed by using pressure values  between subdomains $\Omega_i$ and $\Omega_j$. To be specific, $p_h \in \Omega^*$ is required in the original velocity for constructing in Enhanced Velocity MFEM. However, the post-processed pressure leads to the improved velocity and the visual representation is in Fig. \ref{fig:EVMFEM_postprocess_domains_projection}.  
In case of multiscale setting, it is important to be able to approximate better pressure values nearby the interface. The recovery of velocity computation requires three steps
%There were many ideas in engineering community, for example, to employ an geometric interpolation or other methods. 
\begin{enumerate}
	\item Compute locally $\tilde{p}_h$ from given $\left(p_h, \mathbf{u}_h\right)$
	\item Obtain $s_h$ by using Oswald operator
	\item Compute the velocity at the interface using the two-point flux scheme for $s_h$
\end{enumerate}
We describe construction of $\tilde{p}_h$ and then $s_h$ below.

\subsubsection*{Construction of $\tilde{p}_h$.}
In the Enhance Velocity setting, we may identify $\widehat{\mathbf{V}}_h$ be spaces omitting interface constraints $\mathbf{V}^{\Gamma}$, so $\widehat{\mathbf{V}}_{h,i} :=\bigoplus_{i=1}^n \mathbf{V}_{h,i}(T)$ and then $\widehat{\mathbf{V}}_h :=\bigoplus_{i=1}^n \widehat{\mathbf{V}}_{h,i}$.
Let $\mathbf{u}_h$, $p_h$ be the solution of equations (\ref{eq:EV_eq1})-(\ref{eq:EV_eq2}).
Initially, Lagrange multipliers can be computed in each element. In other words, we define $\lambda_{h, T} \in \Lambda_h$, which is piecewise constant polynomials at edge or face, 
\begin{equation}\label{eq:lagrangemultipliera_fluxreconstr}
\langle \lambda_{h, T}, \mathbf{v}_h \cdot \mathbf{n}_T  \rangle_e :=
\left(\mathbf{K}^{-1}\mathbf{u}_h,\rm  v_h\right)_{T}-\left(p_h, \nabla \cdot \mathbf{v}_h \right)_{T} \qquad \forall \mathbf{v}_h \in  \widehat{\mathbf{V}}_h\left( T \right)
\end{equation}
where the element $T \in \mathcal{T}_{h} $ and its side $e$.
% We note that basis functions are same for $\mathbf{V}^{\Gamma}$ and $\widehat{\mathbf{V}}_h$.
We employ the $L^2$ projected velocity from the interface, which has a finer enhanced velocity approximation, to the edge or face of subdomain element and the formulation is provided in Subsection \ref{sec:diffEV}. We denote polynomial space $\widetilde{W}_h$ in the following manner
\begin{align}
\widetilde{W}_h=\{\varphi_h : \langle \llbracket  \varphi_h \rrbracket, \psi_h \rangle_e=0 \qquad \forall e \in \mathcal{E}^{int}_h \cup \mathcal{E}^{ext}_h  , \forall \psi_h \in \mathbb{Q}_m(e) \}
\end{align}
where $\mathbb{Q}_m$ is standard polynomial space that is defined in \cite{arbogast1995implementation,pencheva2013robust,chen2006computational}.
We next set the post-processed $\tilde{p}_h$ which is proposed in \cite{arbogast1995implementation} and the construction is performed with the following properties, for each $T \in \mathcal{T}_{h} $

\begin{align} 
(\tilde{p}_h, w_h)_T&=(p_h, w_h)_T \qquad \forall w_h \in \widetilde{W}_h(T) \label{eqn:4_32aa} \\
\langle \tilde{p}_h, \mu_h\rangle_e&=\langle \lambda_h, \mu_h \rangle_e \qquad \forall \mu_h \in \Lambda_h(e), \forall e \in \partial T \label{eqn:4_33aa}.	
\end{align}

\subsubsection*{Construction of $s_h$.}
We propose to construct $s_h$ in each subdomain $\Omega_i$ that has the conforming mesh in order to be a computational efficient. Construction of $s_h$  involves the averaging operator $\mathcal{I}_{\rm av} :\mathbb{Q}_k(\mathcal{T}_h) \rightarrow \mathbb{Q}_k(\mathcal{T}_h) \cap H_0^1(\Omega_i) $. For definition of $\mathbb{Q}_m$ we refer reader to \cite{chen2006computational}. The operator is called Oswald operator and appeared in \cite{pencheva2013robust,ainsworth2005robust,karakashian2003posteriori,vohralik2010unified} and the analysis can be found in \cite{burman2007continuous,karakashian2003posteriori}. It is interesting to note that the mapping of the gradient of pressure through Oswald operator also considered in \cite{zienkiewicz1987simple}.
For given $\varphi_h \in \mathbb{Q}_m(\mathcal{T}_{h})$, we regard  the values of $\mathcal{I}_{\rm av}(\varphi_h)$ as being defined at a Lagrange node $V \in \Omega$ by averaging $\varphi_h$ values associated this node,  
\begin{align} \label{eq:oswaldavg}
\mathcal{I}_{\rm av}(\varphi_h)(V)=\frac{1}{\vert \mathcal{T}_{h} \vert} \sum_{T \in \mathcal{T}_{h}} \varphi_h \vert_T (V)
\end{align}
where $\vert  A \vert$ is cardinality of sets $A$ and $\mathcal{T}_{h} $ is all collection of $T \in \mathcal{T}_h$ for fixed $V$. One can see that $\mathcal{T}_{h}(V)=\varphi(V)$ at those nodes that are inside of given $T \in \mathcal{T}_{h}$. We set the value of $\mathcal{I}_{\rm av}(\varphi_h)$ is zero at boundary nodes. 
Now in our setting we define recovered pressure $s_h$ for the locally post-processed $\tilde{p}_h$ as follows.
\begin{equation*}
s_h:=\mathcal{I}_{\rm av}(\tilde{p}_h)
\end{equation*}

\subsection{Implementation steps}
For simplicity, we provide key steps of numerical implementation of post-processed pressure in two dimensional case. However, it can be extended for general cases. Based on piecewise pressure and velocity from the lowest order Raviart-Thomas spaces over rectangles our aim to reconstruct smoother pressure $s_h$.
For given element $T \in \mathcal{T}_h(\Omega_i)$, the main steps are
\begin{enumerate}
	\item Evaluate $\lambda_{h, T}$ at edge $e_j$, $j=1,..4$ based on $(\mathbf{u}_h, p_h)$,
	\item Compute $\tilde{p}_h$ from known $\lambda_{h, T}$, and  $p_h$ by using (\ref{eq:lagrangemultipliera_fluxreconstr}),
	\item Based on $\tilde{p}_h$ compute $s_h$ equation (\ref{eq:oswaldavg}) at Lagrange nodes in $\Omega_i$.
\end{enumerate}

Step 1 is standard computation of Lagrange multiplier for each element. In step 2, we are relying on higher order polynomial, in our case, it is Span$\{1, x, y, x^2, y^2\}$. It is sufficient to store coefficients of polynomials. In step 3, we use \\ Span$\{1, x, y, x^2, y^2, xy, x^2y, xy^2, x^2y^2\}$ and 9 Lagrange nodes of rectangle elements that are four rectangle nodes, four midpoints at edge and center of rectangle. This case each node requires to find neighboring elements values to compute coefficients of $s_h$.

\section{Numerical Examples}\label{numericalexample}
\begin{figure}[H]
	\centering
	\includegraphics[width=0.5\linewidth]{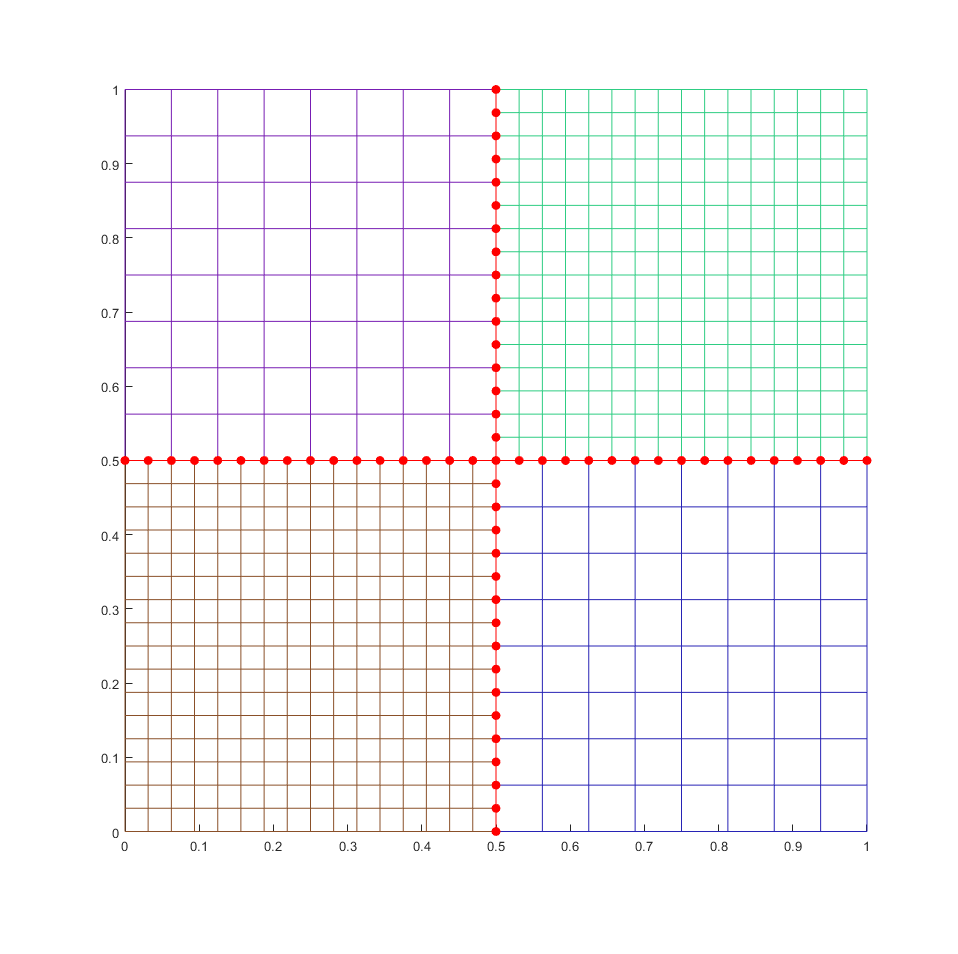}
	\caption{Example of non-matching grids for subdomains.}
	\label{fig:posteriorisubdomainmesh}
\end{figure}

In this section, numerical results are presented to demonstrate challenging problems of velocity approximation at the interface of non-matching multiblock grids.
%We focus on improvement of velocity values at interface by smoothing the pressure values that are adjacent to subdomain interface. 
We have conducted tests for several examples and we concentrate our attention on the interface error for homogeneous and heterogeneous permeability coefficients, respectively. 
%In third and fourth numerical examples, we consider the entire domain error for homogeneous and heterogeneous permeability coefficients, respectively. 
%with diagonal homogeneous permeability coefficients in the numerical tests 1,3  and heterogeneous permeability coefficients in the numerical tests 2,4.
We set same domain $\Omega = \left(0, 1\right) \times \left(0, 1\right)$  for all tests and for some the ratio is $H/h =4$.
% and for others the ratio is $H/h =2$. 
Initial subdomains grids $\mathcal{T}_h$ are chosen in way that has a checkerboard pattern for subdomains. Example of such discretization is shown in Fig. \ref{fig:posteriorisubdomainmesh}. The discrete $L^2$ velocity error 
%$e_{\mathbf{u}_h, \Omega}$ or 
$e_{\mathbf{u}_h, \Gamma}$ is based on the values of the normal component at the midpoint of the edges and is normalized by the analytical solution.

\subsection*{Numerical test 1}
First example tests for uniform permeability, so  $\mathbf{K} = \mathbf{I}$. We report the velocity error and the improved velocity error. We compute the source term and boundary conditions according to the analytical solution, which is taken as follows
\begin{align*}
p(x,y)=\sin(2\pi x)\sin(2\pi y) 
\end{align*}
We set the ratio $H/h=4$ for the result that is shown in below table.
\begin{table}[H] 
	\center
	\begin{tabular}{|l|r|r|r|r|}
		\hline
		$n$ & \multicolumn{2}{ c |}{$e_{\mathbf{u}_h, \Gamma}$}          &  \multicolumn{2}{ c |}{$e_{\mathbf{\tilde{u}}_h, \Gamma}$}      \\
		
		& error & order  & error & order \\
		\hline
		8 & 1.47e-01 & -----   & 3.55e-01   & -----  \\
		\hline
		16 & 7.70e-02  & 0.93 & 1.12e-01  & 1.67 \\
		\hline
		32 & 3.94e-02  & 0.97 & 3.73e-02  & 1.58 \\
		\hline
		48 & 2.65e-02  & 0.98  &  2.09e-02 & 1.43\\
		\hline
	\end{tabular} 
	\vspace{1em}
	\caption{Convergence test 1, velocity and recovered velocity error using the post-processed pressure at interface.}
	\label{table:test1}
\end{table}
%\begin{figure}[!htp]
%	\centering
%	\includegraphics[width=0.7\linewidth]{images/VelocityInterfaceEV/Error_Interface_u_EV_vs_utilde_EV_refin4_test1}
%	\caption{Test1: comparison of flux error at the interface between subdomains to post-processed flux error }
%	\label{fig:error_test1}
%\end{figure}
We list the error of velocity and recovered velocity in Tabel \ref{table:test1}. The results shows the convergence rate improvement for recovered velocity approximation compare to the provided velocity approximation at the interface.
We note that the improvement of convergence rate is significant from order $O(h^{0.95})$ to $O(h^{1.5})$. 
\subsection*{Numerical test 2}
We consider the a diagonal oscillating tensor coefficient as follows.
\begin{align*}
\mathbf{K} =
\begin{bmatrix} 
15-10 \sin(3\pi x)\sin(3 \pi y)    &   0 \\ 
0 &    15-10 \sin(3\pi x)\sin(3 \pi y)
\end{bmatrix} 
\end{align*}
We impose the source term $f$ and Dirichlet boundary condition according to the analytical solution
\begin{align*}
p(x,y) = \sin(2\pi x) \sin(2\pi y).
\end{align*}
We set the ratio $H/h=4$ for the result that is shown in Table \ref{table:test2}.
\begin{table}[H]
	\center
	\begin{tabular}{|l|r|r|r|r|}
		\hline
		
		$n$ & \multicolumn{2}{ c |}{$e_{\mathbf{u}_h, \Gamma}$}          &  \multicolumn{2}{ c |}{$e_{\mathbf{\tilde{u}}_h, \Gamma}$}      \\
		
		& error & order  & error & order \\
		\hline
		8 & 1.78e-01  & ----- &  3.78e-01  & -----\\
		\hline
		16 & 8.89e-02 & 1.00   & 1.00e-01  & 1.91 \\
		\hline
		32 & 4.43e-02 & 1.00  & 2.87e-02   & 1.81 \\
		\hline
		48 & 2.96e-02  & 1.00  & 1.56e-02  & 1.51 \\
		\hline
	\end{tabular}
	\vspace{1em}
	\caption{Convergence test 2: velocity and recovered velocity error using the post-processed pressure at interface.}
	\label{table:test2}
\end{table}
From Table \ref{table:test2}, we see a significant increase on the convergence rate for recovered velocity while the convergence rate of provided velocity stays $\mathcal{O}(h^{1.0})$. We observe that the numerical method is an effective way to improve velocity at the interface between subdomains.

%----------------------------------------------------------------------------------------------------------------------------------
\section{Conclusion}\label{conclusion}
The present study of velocity in Enhanced Velocity Mixed Finite Element Method was designed to investigate the effect of the post-processed pressure on velocity in the interface of subdomains. In this paper, the focus of attention is on the incompressible Darcy flow in the non-matching multiblock grid setting. Multiple numerical results demonstrate that the interface velocity approximation can be improved with using the post-processed pressure. These findings can contribute in several ways to our approximation of velocity and provide a good construction of velocity for a \textit{posteriori} error analysis such as the recovery-based estimate.

%----------------------------------------------------------------------------------------------------------------------------------
\section{Acknowledgments}
\label{sect:acks}

First author would like to thank Drs. I. Yotov and T. Arbogast for discussions on formulation of the different view of EVMFEM. This research is supported by Faculty Development Grant, Nazarbayev University.

%We are thankful for discussion with Professor Todd Arbogast. 
\bibliographystyle{unsrt}  %plainnat 
%\bibliographystyle{splncs04}
%\bibliography{recovery_flux_manuscript}
\bibliography{ref}

\begin{thebibliography}{10}

\bibitem{amanbek2017adaptive}
Yerlan Amanbek, Gurpreet Singh, Mary~F Wheeler, and Hans van Duijn.
\newblock Adaptive numerical homogenization for upscaling single phase flow and
  transport.
\newblock {\em ICES Report}, 12:17, 2017.

\bibitem{singh2017adaptive}
Gurpreet Singh, Yerlan Amanbek, and Mary~F Wheeler.
\newblock Adaptive homogenization for upscaling heterogeneous porous medium.
\newblock In {\em SPE Annual Technical Conference and Exhibition}. Society of
  Petroleum Engineers, 2017.

\bibitem{amanbek2018priori}
Yerlan Amanbek and Mary Wheeler.
\newblock A priori error analysis for transient problems using enhanced
  velocity approach in the discrete-time setting.
\newblock {\em arXiv preprint arXiv:1812.04809}, 2018.

\bibitem{amanbek2018new}
Yerlan Amanbek.
\newblock {\em A new adaptive modeling of flow and transport in porous media
  using an enhanced velocity scheme}.
\newblock PhD thesis, 2018.

\bibitem{gerritsen2008integration}
M~Gerritsen and JV~Lambers.
\newblock Integration of local--global upscaling and grid adaptivity for
  simulation of subsurface flow in heterogeneous formations.
\newblock {\em Computational Geosciences}, 12(2):193--208, 2008.

\bibitem{arbogast2014posteriori}
Todd Arbogast, D~Estep, B~Sheehan, and S~Tavener.
\newblock A posteriori error estimates for mixed finite element and finite
  volume methods for problems coupled through a boundary with nonmatching
  grids.
\newblock {\em IMA Journal of Numerical Analysis}, 34(4):1625--1653, 2014.

\bibitem{wheeler2002enhanced}
John~A Wheeler, Mary~F Wheeler, and Ivan Yotov.
\newblock Enhanced velocity mixed finite element methods for flow in multiblock
  domains.
\newblock {\em Computational Geosciences}, 6(3-4).

\bibitem{thomas2011enhanced}
Sunil~G Thomas and Mary~F Wheeler.
\newblock Enhanced velocity mixed finite element methods for modeling coupled
  flow and transport on non-matching multiblock grids.
\newblock {\em Computational Geosciences}, 15(4):605--625, 2011.

\bibitem{arbogast1995implementation}
Todd Arbogast and Zhangxin Chen.
\newblock On the implementation of mixed methods as nonconforming methods for
  second-order elliptic problems.
\newblock {\em Mathematics of Computation}, 64(211):943--972, 1995.

\bibitem{russell1983finite}
Thomas~F Russell and Mary~Fanett Wheeler.
\newblock {\em Finite element and finite difference methods for continuous
  flows in porous media}, pages 35--106.
\newblock SIAM, 1983.

\bibitem{glowinski1988domain}
Roland Glowinski and Mary~F Wheeler.
\newblock Domain decomposition and mixed finite element methods for elliptic
  problems.
\newblock In {\em First international symposium on domain decomposition methods
  for partial differential equations}, pages 144--172, 1988.

\bibitem{pencheva2013robust}
Gergina~V Pencheva, Martin Vohral{\'\i}k, Mary~F Wheeler, and Tim Wildey.
\newblock Robust a posteriori error control and adaptivity for multiscale,
  multinumerics, and mortar coupling.
\newblock {\em SIAM Journal on Numerical Analysis}, 51(1):526--554, 2013.

\bibitem{ainsworth2005robust}
Mark Ainsworth.
\newblock Robust a posteriori error estimation for nonconforming finite element
  approximation.
\newblock {\em SIAM Journal on Numerical Analysis}, 42(6):2320--2341, 2005.

\bibitem{karakashian2003posteriori}
Ohannes~A Karakashian and Frederic Pascal.
\newblock A posteriori error estimates for a discontinuous galerkin
  approximation of second-order elliptic problems.
\newblock {\em SIAM Journal on Numerical Analysis}, 41(6):2374--2399, 2003.

\bibitem{vohralik2010unified}
Martin Vohral{\'\i}k.
\newblock Unified primal formulation-based a priori and a posteriori error
  analysis of mixed finite element methods.
\newblock {\em Mathematics of Computation}, 79(272):2001--2032, 2010.

\bibitem{chen2006computational}
Zhangxin Chen, Guanren Huan, and Yuanle Ma.
\newblock {\em Computational methods for multiphase flows in porous media},
  volume~2.
\newblock Siam, 2006.

\bibitem{burman2007continuous}
Erik Burman and Alexandre Ern.
\newblock Continuous interior penalty hp-finite element methods for advection
  and advection-diffusion equations.
\newblock {\em Mathematics of Computation}, 76(259):1119--1140, 2007.

\bibitem{zienkiewicz1987simple}
Olgierd~C Zienkiewicz and Jian~Z Zhu.
\newblock A simple error estimator and adaptive procedure for practical
  engineerng analysis.
\newblock {\em International Journal for Numerical Methods in Engineering},
  24(2):337--357, 1987.

\end{thebibliography}

\end{document}